\documentclass{amsart}
\usepackage{graphics}

\newtheorem{theorem}{Theorem}[section]
\newtheorem{lemma}[theorem]{Lemma}
\newtheorem{corollary}[theorem]{Corollary}

\theoremstyle{definition}

\theoremstyle{remark}

\numberwithin{equation}{section}

\DeclareSymbolFont{AMSb}{U}{msb}{m}{n}
\DeclareMathSymbol{\Z}{\mathalpha}{AMSb}{"5A}
\DeclareMathSymbol{\Q}{\mathalpha}{AMSb}{"51}

\begin{document}
\newcommand{\beqs}{\begin{equation*}}
\newcommand{\eeqs}{\end{equation*}}
\newcommand{\beq}{\begin{equation}}
\newcommand{\eeq}{\end{equation}}
\newcommand\nutwid{\overset {\text{\lower 3pt\hbox{$\sim$}}}\nu}
\newcommand\Mtwid{\overset {\text{\lower 3pt\hbox{$\sim$}}}M}
\newcommand\ptwid{\overset {\text{\lower 3pt\hbox{$\sim$}}}p}
\newcommand\pitwid{\overset {\text{\lower 3pt\hbox{$\sim$}}}\pi}
\newcommand\bijone{\overset {1}\longrightarrow}                   
\newcommand\bijtwo{\overset {2}\longrightarrow}                   
\newcommand\pihat{\widehat{\pi}}

\newcommand{\bin}[2]{\genfrac{(}{)}{0pt}{0}{#1}{#2}}

\newcommand\mymod[1]{\mbox{ mod}\ {#1}}
\newcommand\myto{\to}
\newcommand\srank{\mathrm{srank}}
\newcommand\pitc{\pi_{\mbox{\scriptsize $t$-core}}}
\newcommand\pifourc{\pi_{{\mbox{\scriptsize $4$-core}}}}
\newcommand\pitcb[1]{\pi_{\mbox{\scriptsize ${#1}$-core}}}  
\newcommand\mpitcb[1]{\pi_{\mathrm{\scriptsize {#1}-core}}}  
\newcommand\lamseq{(\lambda_1, \lambda_2, \dots, \lambda_\nu)}
\newcommand\mylabel[1]{\label{#1}}
\newcommand\eqn[1]{(\ref{eq:#1})}
\newcommand\stc{{St-crank}}
\newcommand\mstc{\mbox{St-crank}}
\newcommand\mmstc{\mbox{\scriptsize\rm St-crank}}
\newcommand\tqr{{$2$-quotient-rank}}
\newcommand\mtqr{\mbox{$2$-quotient-rank}}
\newcommand\mmtqr{\mathrm{2-quotient-rank}}
\newcommand\fcc{{$5$-core-crank}} 
\newcommand\fc{{5\mbox{-core}}} 
\newcommand\tc{{t\mbox{-core}}}
\newcommand\pifc{\pi_{\scriptsize 5\mbox{-core}}}
\newcommand\bgr{BG-rank} 
\newcommand\mbgr{\mbox{BG-rank}} 
 \newcommand\gbgr{GBG-rank} 
\newcommand\mgbgr{\mbox{GBG-rank}} 
\newcommand\nvec{(n_0, n_1, \dots, n_{t-1})}
\newcommand\parity{\mbox{par}}
\newcommand\epart{\mbox{ep}}
\newcommand\pbar{\overline{p}}
\newcommand\Pbar{\overline{P}}
\newcommand\abar{\overline{a}}
\newcommand\legendre[2]{\genfrac{(}{)}{}{}{#1}{#2}}

\title[The GBG-rank and $t$-cores I. Counting and $4$-cores]{The GBG-rank and $t$-cores I. \\ Counting and $4$-cores}

\author{Alexander Berkovich}
\address{Department of Mathematics, University of Florida, Gainesville,
Florida 32611-8105}
\email{alexb@math.ufl.edu}          
\thanks{Research of both authors was supported in part by NSA grant MSPF-06G-150.}

\author{Frank G. Garvan}
\address{Department of Mathematics, University of Florida, Gainesville,
Florida 32611-8105}
\email{fgarvan@math.ufl.edu}          


\date{August 13, 2008}  



\begin{abstract}
Let $r_j(\pi,s)$ denote the number of cells, colored $j$, in the $s$-residue diagram of partition $\pi$. 
The $\mgbgr$ of $\pi\mymod{s}$ is defined as 
\beqs
\mgbgr(\pi,s) = \sum_{j=0}^{s-1}r_j(\pi,s)e^{I\frac{2\pi}{s}j}, \quad\quad\quad I=\sqrt -1.
\eeqs
We will prove that for $(s,t)=1$
\beqs
\nu(s,t)\leq \frac{\bin{s+t}{s}}{s+t},
\eeqs
where $\nu(s,t)$ denotes a number of distinct values that 
$\mgbgr$ of a $t$-core mod $s$  may assume. 
The above inequality becomes an equality when $s$ is prime or when $s$ is composite and $t\le 2p_s$, 
where $p_s$ is a smallest prime divisor of $s$.\\
We will show that the generating functions for $4$-cores with 
prescribed GBG-rank mod $3$ value are all
$eta$-products.
\end{abstract}
\maketitle

\section{Introduction} \label{sec:intro}
\bigskip

A partition $\pi$ is a nonincreasing sequence
\beqs
\pi =\lamseq
\eeqs
of positive integers (parts) $\lambda_1\geq\lambda_2\geq\lambda_3\geq\cdots\geq\lambda_\nu>0$.
The norm of $\pi$, denoted $|\pi|$, is defined as
\beqs
|\pi| = \sum_{i=1}^{}\lambda_i.
\eeqs
If $|\pi|=n$, we say that $\pi$ is a partition of $n$. The (Young) diagram of $\pi$
is a convenient way to represent $\pi$ graphically: the parts of $\pi$ are shown as rows
of unit squares (cells). Given the diagram of $\pi$ we label a cell in the $i$-th row
and $j$-th column by the least nonnegative integer $\equiv j-i\mymod{s}$. The resulting
diagram is called an $s$-residue diagram \cite{JK}. One can also label cells in the infinite column 
$0$ and the infinite row $0$ in the same fashion. The resulting diagram is called the 
extended $s$-residue diagram of $\pi$ \cite{GKS}.
With each $\pi$ we can associate the $s$-dimensional vector
\beqs
\mathbf r (\pi,s)=(r_0,r_1,\ldots,r_{s-1}),
\eeqs
where $r_i$, $0\leq i\leq s-1$ is the number of cells of $\pi$ labelled $i$ in the $s$-residue diagram of $\pi$.
We shall also require 
\beqs
\mathbf n (\pi,s)=(n_0,n_1,\ldots,n_{s-1}),
\eeqs
where for $0\leq i\leq s-2$ 
\beqs
n_i=r_i-r_{i+1},
\eeqs
and
\beqs
n_{s-1}=r_s-r_0.
\eeqs
Note that 
\beqs
\mathbf n \cdot \mathbf l_s=\sum_{i=0}^{s-1}n_i=0,
\eeqs
where 
\beqs
\mathbf l_s=(1,1,1,\ldots,1)\in \Z^s.
\eeqs
If some cell of $\pi$ shares a vertex or edge with the rim of the diagram of $\pi$, 
we call this cell a rim cell of $\pi$. A connected collection of rim cells of $\pi$ is
called a rim hook of $\pi$ if $\pi \backslash$(rim hook) is a legitimate partition. 
We say that $\pi$ is a $t$-core, denoted $\pitc$, 
if its diagram has no rim hooks of length $t$ \cite{JK}.
The Durfee square of $\pi$, denoted $D(\pi)$, is the largest square that fits inside the diagram of $\pi$.
Reflecting the diagram of $\pi$ about its main diagonal, one gets the diagram of $\pi^\ast$
(the conjugate of $\pi$). More formally, 
\beqs
\pi^\ast = (\lambda_1^\ast,\lambda_2^\ast,\ldots)
\eeqs
with $\lambda_i^\ast$ being the number of parts of $\pi\geq i$.
Clearly,
\beqs
D(\pi)=D(\pi^\ast).
\eeqs
In \cite{BG} we defined a new partition statistic of $\pi$
\beq
\mgbgr(\pi,s): = \sum_{i=0}^{s-1}r_i(\pi,s)\omega_s^i,
\mylabel{eq:1.1}
\eeq
where
\beqs
\omega_s=e^{I\frac{2\pi}{s}}
\eeqs
and
\beqs
I=\sqrt{-1}.
\eeqs
We refer to this statistic as the GBG-rank of $\pi$ mod $s$.
The special case $s=2$ was studied in great detail in \cite{BG} and \cite{BY}. 
In particular, we have shown in \cite{BG} that for any odd $t>1$
\beq
\mgbgr(\pitc,2) = \frac{1-\sum_{i=0}^{t-1}(-1)^{i+n_i(\pitc,t)}}{4}
\mylabel{eq:1.2}
\eeq
and that
\beq
-\bigg\lfloor\frac{t-1}{4}\bigg\rfloor\leq\mgbgr(\pitc,2)\leq\bigg\lfloor\frac{t+1}{4}\bigg\rfloor,
\mylabel{eq:1.3}
\eeq
where $\lfloor x\rfloor$ is the integer part of $x$.
Our main object here is to prove the following generalizations of (\ref{eq:1.2}) and (\ref{eq:1.3}).

\begin{theorem}\label{t1}
Let $t,s\in\Z_{>1}$ and $(t,s)=1$. Then
\beq
\mgbgr(\pitc,s) =
\frac{\sum_{i=0}^{t-1}\omega_s^{i+1}(\omega_s^{{tn_i(\pitc,t)}}-1)}{(1-\omega_s)(1-\omega_s^t)}
\mylabel{eq:1.4}
\eeq
\end{theorem}

\begin{theorem}\label{t2}
Let $\nu(s,t)$ denote the number of distinct values that $\mgbgr$ of $\pitc\mymod{s}$ may assume. 
Then
\beq
\nu(s,t)\leq \frac{\bin{t+s}{t}}{t+s},
\mylabel{eq:1.5}
\eeq
provided $(s,t)=1$.
\end{theorem}

\begin{theorem}\label{t3}
Let $\nu(s,t)$ be as in Theorem 1.2 and $(s,t)=1$. Then
\beq
\nu(s,t) = \frac{\bin{t+s}{t}}{t+s},
\mylabel{eq:1.6}
\eeq
iff either $s$ is prime or $s$ is composite and $t<2p_s$, where $p_s$ is a smallest prime divisor of $s$.
\end{theorem}
\noindent
Our of proof of this Theorem depends crucially on the following

\begin{lemma}\label{l1}
Let $s,t\in\Z_{>1}$ and $(s,t)=1$.
Let $\mathbf j=(j_0,j_1,\ldots,j_{t-1})$, $\mathbf{\tilde j}=(\tilde j_0, \tilde j_1,\ldots,\tilde j_{t-1})$ be integer 
valued vectors such that
\beq
0\leq j_0\leq j_1\leq\cdots\leq j_{t-1}\leq s-1,
\mylabel{eq:1.7}
\eeq
\beq
0\leq\tilde j_0\leq\tilde j_1\leq\cdots\leq\tilde j_{t-1}\leq s-1,
\mylabel{eq:1.8}
\eeq
and
\beq
\sum_{i=0}^{t-1} \omega_s^{j_i} = \sum_{i=0}^{t-1} \omega_s^{\tilde j_i} 
\mylabel{eq:1.9}
\eeq
\beq
\prod_{i=0}^{t-1} \omega_s^{j_i} = \prod_{i=0}^{t-1} \omega_s^{\tilde j_i}.
\mylabel{eq:1.10}
\eeq
Then
\beqs
\mathbf j = \mathbf{\tilde j},
\eeqs
iff either $s$ is prime or $s$ is composite such that $t<2p_s$, where $p_s$ is
a smallest prime divisor of $s$.
\end{lemma}

The rest of this paper is organised as follows. In Section 2, we collect some necessary background on $t$-cores 
and prove Theorems 1.1 and Theorem 1.2. Section 1.3 is devoted to the proof of Lemma 1.4 and Theorem 1.3. 
Section 4 deals with $4$-cores with prescribed values of $\mgbgr\mymod{3}$. 
There we will provide new combinatorial interpretation and proof of the 
Hirshhorn-Sellers identities for $4$-cores \cite{HS}.
We conclude with the remarks connecting this development and that of \cite{O} and \cite{A}.

\bigskip
\section{Properties of the GBG-rank} \label{sec:2}
\bigskip

We begin with some definitions from \cite{GKS}.
A region $r$ in the extended $t$-residue diagram of $\pi$ is the set of all cells $(i,j)$ satisfying $t(r-1)\leq j-i<tr$. 
A cell of $\pi$ is called exposed if it is at the end of a row of $\pi$. One can construct $t$ bi-infinite words 
$W_0,W_1,\ldots,W_{t-1}$ of two letters $N,E$ as follows:
The $r$th letter of $W_i$ is $E$ if there is an exposed cell labelled $i$ in the region $r$ of $\pi$, 
otherwise the $r$th letter of $W_i$ is $N$. It is easy to see that the word set $\{W_0,W_1,\ldots,W_{t-1}\}$ 
fixes $\pi$ uniquely. It was shown in \cite{GKS} that $\pi$ is a $t$-core iff each word of $\pi$ is of the form:
\begin{align}
\mbox{Region}&:\quad\cdots\cdots\;n_{i-1}\quad n_i\quad n_{i+1}\quad n_{i+2}\quad\cdots\cdots \nonumber\\
          W_0&:\quad\cdots\cdots\quad E \quad\; E \quad\quad N \quad\quad N \;\quad\cdots\cdots.
\mylabel{eq:2.1}
\end{align}
For example, the word image of $\pi_{3\mbox{-core}}=(4,2)$ is 
\begin{align*}
\mbox{Region}&:\quad\cdots\cdots\;\;-1\quad 0\quad 1\quad\;\; 2\quad\;\; 3\quad\cdots\cdots \\
          W_0&:\quad\cdots\cdots\quad E \quad E \quad E \quad E \quad N \quad\cdots\cdots \\
          W_1&:\quad\cdots\cdots\quad E \quad N \quad N \quad N \quad N \quad\cdots\cdots \\
          W_2&:\quad\cdots\cdots\quad E \quad N \quad N \quad N \quad N \quad\cdots\cdots,
\end{align*}
while the associated $\mathbf r$ and $\mathbf n$ vectors are 
$\mathbf r =(r_0,r_1,r_2)=(3,1,2)$, $\mathbf n =(n_0,n_1,n_2)=(2,-1,-1)$, respectively.
In general, the map
$$
\phi(\pitc)=\mathbf n(\pitc,t) =(n_0,n_1,\dots,n_{t-1})
$$
is a bijection 
from the set of $t$-cores to the set
$$
\{\mathbf n \in \Z^t\,:\, \mathbf n \cdot \mathbf 1 = 0\}.
$$
Next, we mention three more useful facts from \cite{GKS}.\\ 
\noindent
\textsl{A.}
\beq
\big|\pitc\big|=\frac{t}{2}|\mathbf n |^2+\mathbf b_t\cdot\mathbf n,
\mylabel{eq:2.2} 
\eeq
where $\mathbf b_t=(0,1,2,\ldots,t-1)$. \\
\noindent
\textsl{B.}
\beq
\sum_{i\in P_1}n_i=-\sum_{i\in P_{-1}}n_i=D(\pitc),
\mylabel{eq:2.3}
\eeq
where $P_\alpha=\{i\in\Z:0\leq i\leq t-1,\alpha n_i>0\}$, $\alpha=-1,1$.\\
\noindent
\textsl{C.}\\
Under conjugation $\phi(\pitc)$ transforms as 
\beq
(n_0,n_1,n_2,\ldots,n_{t-1})\rightarrow(-n_{t-1},-n_{t-2},\ldots,-n_0).   
\mylabel{eq:2.4}
\eeq

We begin our proof of the Theorem 1.1 by observing that under conjugation $\mgbgr$ transforms as
\beq
\mgbgr(\pi,s)=\sum_{i=0}^{s-1} r_i\omega_s^i \Longrightarrow \mgbgr(\pi^\ast,s) = \sum_{i=0}^{s-1} r_i\omega_s^{-i}. 
\mylabel{eq:2.5}
\eeq
Next, we use that
\beq
\mgbgr(\pi,s)= \mgbgr(\pi_1,s)+ \mgbgr(\pi_2,s)- D.
\mylabel{eq:2.6}
\eeq
Here, $\pi_1$ is obtained from the diagram of $\pitc$ by removing all cells strictly below the main diagonal of $\pitc$.
Similarly, $\pi_2$ is obtained from $\pitc$ by removing the cells strictly to the right of the main diagonal.\\
Recalling \eqn{2.1} and \eqn{2.3} we find that 
\beq
\mgbgr(\pi_1,s)=\sum_{i\in P_1}\sum_{k=1}^{n_i}\sum_{j=0}^{i+t(k-1)} \omega_s^j=
\frac{D}{1-\omega_s}-
\sum_{i\in P_1}\frac{\omega_s^{i+1}(1-\omega_s^{tn_i})}{(1-\omega_s)(1-\omega_s^t)}. 
\mylabel{eq:2.7}
\eeq
Analogously,
\beq
\mgbgr(\pi_2^\ast,s)=\frac{D}{1-\omega_s}-\sum_{i\in P_{-1}}\frac{\omega_s^{t-i}(1-\omega_s^{-tn_i})}{(1-\omega_s)(1-\omega_s^t)}, 
\mylabel{eq:2.8}
\eeq
where we made use of \eqn{2.4}.\\
Clearly, \eqn{2.5} and \eqn{2.8} imply that  
\beq
\mgbgr(\pi_2,s)=-\frac{D\omega_s}{1-\omega_s}-
\sum_{i\in P_{-1}}\frac{\omega_s^{1+i}(1-\omega_s^{tn_i})}{(1-\omega_s)(1-\omega_s^t)}. 
\mylabel{eq:2.9}
\eeq
Next, we combine \eqn{2.6}, \eqn{2.7} and \eqn{2.9} to find that  
\beqs
\mgbgr(\pitc,s)=-\sum_{i\in P_{-1}\bigcup P_1} \frac{\omega_s^{1+i}(1-\omega_s^{tn_i})}{(1-\omega_s)(1-\omega_s^t)}
=\sum_{i=0}^{t-1}\frac{\omega_s^{1+i}(\omega_s^{tn_i}-1)}{(1-\omega_s)(1-\omega_s^t)}, 
\eeqs
as desired.\\

Our proof of Theorem 1.2 involves three observations, which we now proceed to discuss.\\
\textsl{Observation 1}:

Let $a_r(s,t)$ denote the number of vectors 
$\mathbf j = (j_0,j_1,\ldots,j_{t-1})$ such that
\beqs
0\leq j_0\leq j_1\leq j_2\leq\cdots\leq j_{t-1}<s,
\eeqs
\beqs
\sum_{k=0}^{t-1}j_k\equiv r \mymod{s}.
\eeqs
Then 
\beq
\nu(s,t)\leq a_{\frac{t(t+1)}{2}}(s,t),
\mylabel{eq:2.10}
\eeq
provided $(s,t)=1$.\\
\noindent
\textbf{Proof.}

Suppose $(s,t)=1$. It is clear that the number of values of the GBG-rank of
$t$-cores mod $s$ is the number distinct values of
$$
\sum_{i=0}^{t-1}\omega_s^{1+i+t n_i},
$$
where $\mathbf n \in \Z^t$ and $\mathbf n \cdot \mathbf 1_t=0$.
Given any such $\mathbf n$-vector we reduce the exponents ${1+i+t n_i}$
mod $s$ and reorder to obtain a $\mathbf j$-vector such that
$$
\sum_{k=0}^{t-1} j_k \equiv \sum_{i=0}^{t-1} {1+i+t n_i}
\equiv \tfrac{t(t+1)}{2}\pmod{s}.
$$
It follows that
$$
\nu(s,t)\leq a_{\frac{t(t+1)}{2}}(s,t).
$$

\noindent
\textsl{Observation 2}:
\beq
\sum_{r=0}^{s-1}a_r(s,t)=\bin{t+s-1}{t}.
\mylabel{eq:2.14}
\eeq
This result is well known and we omit the proof.
Finally, we need \\ 
\noindent
\textsl{Observation 3}:

If $(s,t)=1$ then
\beq
a_0(s,t)=a_1(s,t)=\cdots=a_{s-1}(s,t).
\mylabel{eq:2.15}
\eeq

\noindent
\textbf{Proof.}

There exists an integer $T$ such that $T\cdot t\equiv 1\mymod{s}$, because $s$ and $t$ are coprime. 
This implies that
\beqs
\sum_{i=0}^{t-1}(j_i+T)\equiv 1+\sum_{i=0}^{t-1}j_i\mymod{s}.
\eeqs
Consequently, $a_r(s,t)=a_{r+1}(s,t)$, as desired.
Combining \eqn{2.10}, \eqn{2.14} and \eqn{2.15} we see that 
\beqs
\nu(s,t)\leq a_{\frac{t(t+1)}{2}}=\frac{\bin{s-1+t}{t}}{s}=\frac{\bin{s+t}{t}}{s+t},
\eeqs
and we have Theorem 1.2.

\bigskip
\section{Roots of unity and the number of values of the GBG-rank} \label{sec:3}
\bigskip

It is clear from our proof of Theorem 1.2 that
\beqs
\nu(s,t)=\frac{\bin{s+t}{t}}{s+t}.
\eeqs
iff each $\mathbf j=(j_0,j_1,\ldots,j_{t-1})$ such that
\beqs
0\leq j_0\leq j_1\leq\cdots\leq j_{t-1}<s,
\eeqs
\beqs
\prod_{i=0}^{t-1}\omega_s^{j_i}=\omega_s^{\frac{t(t+1)}{2}}
\eeqs
is associated with a distinct complex number $\sum_{i=0}^{t-1}\omega_s^{j_i}$. Lemma 1.4 tells us when this
is exactly the case. This means that Theorem 1.3 is an immediate corollary of this Lemma. To prove it we need
to consider six cases.\\
\noindent
\textbf{Case 1.}

$s$ is prime, $(s,t)=1$. Note that 
\beqs
\Phi_s(x):=1+x+x^2+\cdots +x^{s-1}
\eeqs
is a minimal polynomial of $\omega_s$ over $\Q$.
Let us now define
\beq
p_1(x):=\sum_{i=0}^{t-1}(x^{j_i}-x^{\tilde j_i}),  
\mylabel{eq:3.1}
\eeq
where $\mathbf j$ and $\mathbf {\tilde j}$ satisfy the constraints \eqn{1.7} - \eqn{1.10}.
It is clear that
\beqs
p_1(\omega_s)=0,
\eeqs
\beqs
p_1(1)=0,
\eeqs
and that $\mbox{deg}(p_1(x))<s$. But $(x-1)\Phi_s(x)$ divides $p_1(x)$. 
This implies that $p_1(x)$ is identically zero and  
$\mathbf j=\mathbf{\tilde j}$, as desired.

\noindent
\textbf{Case 2}.

Here $s$ is composite, $(s,t)=1$ and $t<2p_s$, where $p_s$ is a smallest prime divisor of $s$. 
Once again \eqn{1.9} implies that
\beqs
p_1(\omega_s)=0.
\eeqs
Moreover, the $s$th cyclotonic polynomial, defined as
\beq
\Phi_s(x):=\prod_{\substack{0<j<s,\\(j,s)=1}}(x-\omega_s^j),
\mylabel{eq:3.2}
\eeq
is a minimal polynomial of $\omega_s$ over $\Q$.
This means that
\beq
p_1(\omega_s^m)=0,
\mylabel{eq:3.3}
\eeq
for any $1\leq m<s$ such that $(s,m)=1$.
In particular, we have that 
\beq
p_1(\omega_s^k)=0, \quad\quad\quad 1\leq k\leq p_s-1.
\mylabel{eq:3.4}
\eeq
At this point, it is expedient to rewrite \eqn{3.4} as
\beq
h_k(\omega_s^{j_0},\ldots,\omega_s^{j_{t-1}})=h_k(\omega_s^{\tilde j_0},\ldots,\omega_s^{\tilde j_{t-1}}), 
\quad\quad 1\leq k\leq p_s-1,
\mylabel{eq:3.5}
\eeq
where
\beqs
h_k(x_1,x_2,\ldots,x_t)=x_1^k+x_2^k+\cdots +x_t^k.
\eeqs
Next, we use Newton's theorem on symmetric polynomials to convert \eqn{3.5} into $p_s-1$ identities
\beq
\sigma_k(\omega_s^{j_0},\ldots,\omega_s^{j_{t-1}})=\sigma_k(\omega_s^{\tilde j_0},\ldots,\omega_s^{\tilde j_{t-1}}), 
\quad\quad 1\leq k\leq p_s-1,
\mylabel{eq:3.6}
\eeq
where the $k$th elementary symmetric polynomials $\sigma_k$'s in $x_1,x_2,\ldots,x_t$ are defined in a standard way as
\beq
\sigma_k(x_1,x_2,\ldots,x_t)
=\sum_{1\leq i_1\leq i_2<\ldots<i_k\leq t}x_{i_1} x_{i_2}\cdots x_{i_k},
\quad\quad 1\leq k\leq t.
\mylabel{eq:3.7}
\eeq
Note that we can rewrite \eqn{1.10} now as 
\beq
\sigma_t(\omega_s^{j_0},\ldots,\omega_s^{j_{t-1}})=\sigma_t(\omega_s^{\tilde j_0},\ldots,\omega_s^{\tilde j_{t-1}}). 
\mylabel{eq:3.8}
\eeq
But
\beq
\sigma_t\sigma_k^\ast=\sigma_{t-k}, 
\mylabel{eq:3.9}
\eeq
where
$$
\sigma_k^\ast (x_1,x_2,\ldots,x_t) =
\sigma_k(x_1^{-1},x_2^{-1},\ldots,x_{t}^{-1}).
$$
This fortunate fact enables us to convert \eqn{3.6} into $p_s-1$ identities
\beq
\sigma_k(\omega_s^{j_0},\ldots,\omega_s^{j_{t-1}})=\sigma_k(\omega_s^{\tilde j_0},\ldots,\omega_s^{\tilde j_{t-1}}), 
\quad\quad t-p_s+1\leq k\leq t-1.
\mylabel{eq:3.10}
\eeq
But $t<2p_s$, and so, $t-p_s+1\leq p_s$. This means that we have the following $t$ identities
\beq
\sigma_k(\omega_s^{j_0},\ldots,\omega_s^{j_{t-1}})=\sigma_k(\omega_s^{\tilde j_0},\ldots,\omega_s^{\tilde j_{t-1}}), 
\quad\quad 1\leq k\leq t.
\mylabel{eq:3.11}
\eeq
Consequently,
\beqs
\prod_{i=0}^{t-1}(x-\omega_s^{j_i})=\prod_{i=0}^{t-1}(x-\omega_s^{\tilde j_i}).
\eeqs
Recalling that $\mathbf j,\mathbf{\tilde j}$ satisfy \eqn{1.7} and \eqn{1.8}, 
we conclude that $\mathbf j=\mathbf{\tilde j}$. \\
Let us summarize. If $s$ is a prime or if $s$ is a composite number such that $t<2p_s$, then $\mathbf j=\mathbf{\tilde j}$,
provided that $(s,t)=1$ and $\mathbf j,\mathbf{\tilde j}$ satisfy \eqn{1.7}--\eqn{1.10}.\\
It remains to show that $\mathbf j=\mathbf{\tilde j}$ does not have to be true if $s$ is a composite number 
and $t\geq 2p_s$. To this end consider 
\begin{align*}
\mathbf j: & =(0,0,\ldots,0,1,1,3,3)\in\Z^t, \\
\mathbf{\tilde j}: & =(0,0,\ldots,0,0,0,2,2)\in\Z^t, 
\intertext{if $s=4$, $t\geq 4$,}
\mathbf j: & =(0,0,\ldots,0,1,1,4,4)\in\Z^t, \\
\mathbf{\tilde j}: & =(0,0,\ldots,0,0,2,3,5)\in\Z^t, 
\intertext{if $s=6$, $t\geq 4$ and}
\mathbf j: & =(0,0,\ldots,0,3,3,6,6)\in\Z^t, \\
\mathbf{\tilde j}: & =(0,0,\ldots,0,1,2,4,5,7,8)\in\Z^t, 
\end{align*}
if $s=9$, $t\geq 6$, respectively.
It is not hard to verify that $\mathbf j,\mathbf{\tilde j}$ satisfy \eqn{1.7}--\eqn{1.10} and that 
$\mathbf j\neq\mathbf{\tilde j}$ in these cases.
It remains to consider the last case where $s$ is a composite number $\neq 4,6,9$, $t\geq 2p_s$. In this case $s>3p_s$. 
And, as a result,
\beqs
3+\frac{s}{p_s}(p_s-1)<s.
\eeqs
Let us now consider
\begin{align*}
\mathbf j: & =(0,0,\ldots,0,2,2,2+\frac{s}{p_s},2+\frac{s}{p_s},\ldots,2+\frac{s}{p_s}(p_s-1),2+\frac{s}{p_s}(p_s-1))
\in\Z^t, \\
\mathbf{\tilde j}: & =(0,0,\ldots,0,1,3,1+\frac{s}{p_s},3+\frac{s}{p_s},\ldots,1+\frac{s}{p_s}(p_s-1),3+\frac{s}{p_s}(p_s-1))
\in\Z^t. 
\end{align*}
Again, it is straightforward to check that $\mathbf j,\mathbf{\tilde j}$ satisfy 
\eqn{1.7}--\eqn{1.10} and that $\mathbf j\neq\mathbf{\tilde j}$.
This completes our proof of Lemma 1.4.

We have an immediate 
\begin{corollary} Theorem 1.3 holds true.
\end{corollary}
To illustrate the usefulness of Theorem 1.3, consider the following example: $s=3$, $t=4$.
In this case we should have exactly $\frac{\bin{4+3}{3}}{4+3}=5$ distinct values of $\mgbgr(\pi_{4\mbox{-core}},3)$.
To determine these distinct values we substitute the following $\mathbf n$-vectors 
$(0,-1,1,0)$, $(0,0,0,0)$, $(-1,0,0,1)$, $(0,0,-1,1)$, $(-1,1,0,0)$ into \eqn{1.4} 
to obtain $-1,0,1,-\omega_3,-\omega_3^2$, respectively.
To verify this we note that there are exactly $27$ vectors such that
\beqs
\mathbf n\in\Z_3^4 \quad \mbox{and} \quad \mathbf n\cdot\mathbf 1_4\equiv 0\mymod{3}.
\eeqs
In Table 1 we list all these vectors together with the associated $\mgbgr\mymod{3}$ values, determined by \eqn{1.4}. 
These vectors will come in handy later.
\begin{table}
\begin{tabular}{|l|c|}
\hline
$\quad\quad\mathbf n$ vectors & $\mgbgr$ values \\
\hline
$\mathbf n_1=(0,-1,1,0)$ &  -1 \\ 
\hline
$\mathbf n_2=(0,0,0,0)$ & 0 \\ 
$\mathbf n_3=(1,1,-2,0)$ & 0 \\ 
$\mathbf n_4=(-1,-1,1,1)$ & 0 \\ 
$\mathbf n_5=(0,-1,-1,2)$ & 0 \\ 
$\mathbf n_6=(1,-1,0,0)$ & 0 \\ 
$\mathbf n_7=(0,1,-2,1)$ & 0 \\ 
$\mathbf n_8=(2,-1,-1,0)$ & 0 \\ 
$\mathbf n_9=(0,0,1,-1)$ & 0 \\ 
$\mathbf n_{10}=(0,1,-1,0)$ & 0 \\ 
$\mathbf n_{11}=(-1,0,1,0)$ & 0 \\ 
$\mathbf n_{12}=(1,-1,1,-1)$ & 0 \\ 
$\mathbf n_{13}=(0,-1,0,1)$ & 0 \\ 
\hline
$\mathbf n_{14}=(1,1,0,-2)$ & 1 \\ 
$\mathbf n_{15}=(-1,1,-1,1)$ & 1 \\ 
$\mathbf n_{16}=(2,0,-1,-1)$ & 1 \\ 
$\mathbf n_{17}=(1,0,0,-1)$ & 1 \\ 
$\mathbf n_{18}=(1,1,-1,-1)$ & 1 \\ 
$\mathbf n_{19}=(-1,0,0,1)$ & 1 \\ 
\hline
$\mathbf n_{20}=(1,0,-1,0)$ & $-\omega_3$ \\ 
$\mathbf n_{21}=(1,0,-2,1)$ & $-\omega_3$ \\ 
$\mathbf n_{22}=(1,-1,-1,1)$ & $-\omega_3$ \\ 
$\mathbf n_{23}=(0,0,-1,1)$ & $-\omega_3$ \\ 
\hline
$\mathbf n_{24}=(-1,1,0,0)$ & $-\omega_3^2$ \\ 
$\mathbf n_{25}=(-1,1,1,-1)$ & $-\omega_3^2$ \\ 
$\mathbf n_{26}=(-1,2,0,-1)$ & $-\omega_3^2$ \\ 
$\mathbf n_{27}=(0,1,0,-1)$ & $-\omega_3^2$ \\ 
\hline
\end{tabular}
\medskip
\caption{}
\label{Table1}
\end{table}

\bigskip
\section{The GBG-rank of $4$-cores mod $3$} \label{sec:4}
\bigskip

Let $G_t(q)$ denote the generating function for $t$-cores.
\beq
G_t(q):=\sum_{\pitc}q^{|\pitc|}.
\mylabel{eq:4.1}
\eeq
Let $P$ be the set of all partitions and $P_{t\mbox{-core}}$ be the set of all $t$-cores.
There is a well-known bijection
\beqs
\tilde\phi:P\rightarrow P_{t\mbox{-core}}\times P \times P \times P \ldots\times P 
\eeqs
which goes back to D.E. Littlewood \cite{L}
\beqs
\tilde\phi(\pi)=(\pi_{t\mbox{-core}},\hat\pi_0,\hat\pi_1,\ldots,\hat\pi_{t-1})
\eeqs
such that
\beqs
|\pi|=|\pi_{t\mbox{-core}}|+t\sum_{i=0}^{t-1}|\hat\pi_i|.
\eeqs
The multipartition $(\hat\pi_0,\hat\pi_1,\ldots,\hat\pi_{t-1})$ is called the $t$-quotient of $\pi$. 
The immediate corollary of the Littlewood bijection is 
\beq
G_t(q)=\frac{E^t(q^t)}{E(q)},
\mylabel{eq:4.2}
\eeq
where
\beq
E(q):=\prod_{j\geq 1}(1-q^j).
\mylabel{eq:4.3}
\eeq
On the other hand, formula \eqn{2.2} suggests \cite{GKS} that
\beq
G_t(q)=\sum_{\substack{\mathbf n\in\Z^t, \\ \mathbf n \cdot\mathbf 1_t=0}}q^
{\frac{t}{2}|\mathbf n|^2 + \mathbf n\cdot\mathbf b_t},
\mylabel{eq:4.4}
\eeq
so that
\beq
\sum_{\substack{\mathbf n\in\Z^t, \\ \mathbf n \cdot\mathbf 1_t=0}}q^
{\frac{t}{2}|\mathbf n|^2 + \mathbf n\cdot\mathbf b_t}=\frac{E^t(q^t)}{E(q)}.
\mylabel{eq:4.5}
\eeq
The above identity was first obtained by Klyachko \cite{Kl}, who observed that it is a special case of $A_{t-1}$
MacDonald's identity. An elementary proof of \eqn{4.5} can be found in \cite{BG}.
Next we define
\beq
g_c(q)=\sum_{\substack{\pifourc,\\ \mgbgr(\pifourc,3)=c}}q^{|\pifourc|}
\mylabel{eq:4.6}
\eeq
In other words, $g_c(q)$ is the generating function for $4$-cores with a given value $c$ of the $\mgbgr\mymod{3}$.
>From the discussion at the end of the last section it is clear that
\beq
\frac{E^4(q^4)}{E(q)}=g_{-1}(q)+g_0(q)+g_1(q)+g_{-\omega_3}(q)+g_{-\omega_3^2}(q).
\mylabel{eq:4.7}
\eeq
It turns out that 
\begin{align}
g_{-1}(q) & = q^5\frac{E^4(q^{36})}{E(q^9)}, \mylabel{eq:4.8} \\
g_0(q) & = \frac{E^6(q^6) E^2(q^{18})}{E^3(q^3) E(q^{12}) E(q^{36})}, \mylabel{eq:4.9} \\
g_1(q) & = q\frac{E^2(q^9) E^4(q^{12})}{E(q^3) E(q^6) E(q^{18})}, \mylabel{eq:4.10} \\
g_{-\omega_3}(q) & = q^2\frac{E^2(q^9) E(q^{12})E(q^{36})}{E(q^3)}, \mylabel{eq:4.11} \\
g_{-\omega_3^2}(q) & = q^2\frac{E^2(q^9) E(q^{12})E(q^{36})}{E(q^3)}. \mylabel{eq:4.12}
\end{align}
Hence  
\beq
\begin{split}
\frac{E^4(q^4)}{E(q)} & =
\frac{E^6(q^6) E^2(q^{18})}{E^3(q^3) E(q^{12}) E(q^{36})} \\
& + q\frac{E^2(q^9) E^4(q^{12})}{E(q^3) E(q^6) E(q^{18})} \\
& + 2q^2\frac{E^2(q^9) E(q^{12})E(q^{36})}{E(q^3)} \\
& + q^5\frac{E^4(q^{36})}{E(q^9)}.
\end{split}
\mylabel{eq:4.13} 
\eeq
We note that the identity
$$
g_{-\omega_3}(q) = g_{-\omega_3^2}(q)
$$
follows from \eqn{2.5} and the fact that $\pi$ is a $t$-core if and only if
the conjugate $\pi^*$ is.
The identities equivalent to \eqn{4.13} were first proven by Hirschhorn and Sellers \cite{HS}.
However, combinatorial identities \eqn{4.7}--\eqn{4.12} given here are brand new.
The proof of \eqn{4.8} is rather simple. Indeed, data in Table 1, suggests that 
\beq
\begin{split}
g_{-1}(q) & =
\sum_{\substack{\mathbf n\cdot\mathbf 1_4=0, \\ \mathbf n\equiv\mathbf n_1\mymod{\Z_3^4}}}
q^{2|\mathbf n|^2+\mathbf b_4\cdot\mathbf n} \\
& = q^5 \sum_{\mathbf{\tilde n}\cdot\mathbf 1_4=\tilde n_0+\tilde n_1+\tilde n_2+\tilde n_3=0}
q^{9(2|\mathbf{\tilde n}|^2-\tilde n_1+2\tilde n_2+\tilde n_3)} \\
& = q^5 \sum_{\mathbf{\tilde n}\cdot\mathbf 1_4=0}
q^{9(2|\mathbf{\tilde n}|^2+\tilde n_0+2\tilde n_3+3\tilde n_2)} \\
& = q^5 \frac{E^4(q^{36})}{E(q^9)}
\end{split}
\mylabel{eq:4.14} 
\eeq
where in the last step we relabelled the summation variables and used \eqn{4.5} with $t=4$ and $q\rightarrow q^9$.

In what follows we shall require the Jacobi triple product identity (\cite{GR},II.28)
\beq
\sum_{n=-\infty}^\infty(-1)^n q^{n^2}z^n = E(q^2)[zq;q^2]_\infty,
\mylabel{eq:4.15} 
\eeq
where
\beqs
[z;q]_\infty:=\prod_{j=0}^\infty(1-zq^j)\bigg(1-\frac{q^{1+j}}{z}\bigg)
\eeqs
and the formula (\cite{GR},ex.5.21)
\beq
\bigg[ux,\frac{u}{x},vy,\frac{v}{y};q\bigg]_\infty=\bigg[uy,\frac{u}{y},vx,\frac{v}{x};q\bigg]_\infty+
\frac{v}{x}\bigg[xy,\frac{x}{y},uv,\frac{u}{v};q\bigg]_\infty,
\mylabel{eq:4.16} 
\eeq
where
\beqs
[z_1,z_2,\ldots,z_n;q]_\infty:=\prod_{j=1}^n[z_i;q]_\infty.
\eeqs
Setting $u=q^5$, $v=q^3$, $x=q^2$, $y=q$ and replacing $q$ by $q^{12}$ in \eqn{4.16} we find that
\beq
[q^2,q^3;q^{12}]_\infty([q^5;q^{12}]_\infty-q[q;q^{12}]_\infty)=[q,q^5,q^6;q^{12}]_\infty.
\mylabel{eq:4.17} 
\eeq
Analogously, \eqn{4.16} with $u=q^5$, $v=q^2$, $x=q$, $y=1$ and $q\rightarrow q^{12}$ becomes
\beq
[q^5;q^{12}]_\infty+q[q;q^{12}]_\infty=\frac{[q^2,q^2,q^4,q^6;q^{12}]_\infty}{[q,q^3,q^5;q^{12}]_\infty}.
\mylabel{eq:4.18} 
\eeq
Finally, setting $u=q^6$, $v=q^4$, $x=q^3$, $y=1$ and $q\rightarrow q^{12}$ in \eqn{4.16} yields
\beq
[q^3,q^4;q^{12}]_\infty^2=[q,q^5,q^6,q^6;q^{12}]_\infty+q[q^2,q^3;q^{12}]_\infty^2.
\mylabel{eq:4.19} 
\eeq
Next, we use again Table 1 to rewrite \eqn{4.9} as
\beq
\sum_{j=2}^{13}\sum_{\substack{\mathbf n\cdot\mathbf 1_4=0, \\ \mathbf n\equiv \mathbf n_j\mymod{\Z_3^4}}}
q^{2|\mathbf n|^2 + \mathbf b_4\cdot\mathbf n}=R_1(q),
\mylabel{eq:4.20} 
\eeq 
where
\beq
R_1(q)=\frac{E^6(q^6) E^2(q^{18})}{E^3(q^3) E(q^{12}) E(q^{36})}.
\mylabel{eq:4.21} 
\eeq 
Remarkably, \eqn{4.20} is a constant term in $z$ of the following identity
\beq
\sum_{j=2}^{13}s_j(z,q)=R_1(q)\sum_{n=-\infty}^\infty q^{9\frac{n(n+1)}{2}}z^n,
\mylabel{eq:4.22} 
\eeq
where 
\beq
s_j(z,q):=\sum_{\mathbf n\equiv \mathbf n_j\mymod{\Z_3^4}}
q^{2|\mathbf n|^2 + \mathbf b_4\cdot\mathbf n} z^{\frac{\mathbf n\cdot\mathbf 1_4}{3}},\quad\quad\quad j=1,2,\ldots,27.
\mylabel{eq:4.23} 
\eeq 
Using simple changes of variables,
it is straightforward to check that 
\beq
zq^9s_i(zq^9,q)=s_j(z,q),
\mylabel{eq:4.24} 
\eeq
holds true for the following $(i,j)$ pairs:
$(2,3)$, $(3,4)$, $(4,5)$, $(5,2)$, $(6,7)$, $(7,8)$, $(8,9)$, $(9,6)$, $(10,11)$, $(11,12)$, $(12,13)$, $(13,10)$, and that
\beq
zq^9\sum_{n=-\infty}^\infty q^{9\frac{n(n+1)}{2}}(zq^9)^n=\sum_{n=-\infty}^\infty q^{9\frac{n(n+1)}{2}}z^n.
\mylabel{eq:4.25} 
\eeq
Consequently, both sides of \eqn{4.22} satisfy the same first order functional equation
\beq
zq^9 f(zq^9,q)=f(z,q).
\mylabel{eq:4.26} 
\eeq
Thus to prove \eqn{4.22} it is sufficient to verify it at one nontrivial point, 
say $z=z_0:=-q^{-6}$.
It is not hard to check that
\beq
s_4(z_0,q)=s_8(z_0,q)=s_{11}(z_0,q)=0,
\mylabel{eq:4.27} 
\eeq
and that
\beq
s_3(z_0,q)+s_9(z_0,q)=s_5(z_0,q)+s_{12}(z_0,q)=0.
\mylabel{eq:4.28} 
\eeq
We see that \eqn{4.22} with $z=z_0$ becomes
\beq
s_2(z_0,q)+s_6(z_0,q)+s_7(z_0,q)+s_{10}(z_0,q)+s_{13}(z_0,q)=R_1(q)\sum_{n=-\infty}^\infty q^{9\frac{n(n+1)}{2}}z_0^n.
\mylabel{eq:4.29} 
\eeq
Upon making repeated use of \eqn{4.15} and replacing $q^3$ by $q$ 
we find that \eqn{4.9} is equivalent to
\beq
\begin{split}
[q^4,q^5,q^5,q^6;q^{12}]_\infty & + q[q^2,q^3,q^4;q^{12}]_\infty([q^5;q^{12}]_\infty-q[q;q^{12}]_\infty) \\ 
& + q[q,q^4,q^5,q^6;q^{12}]_\infty + q^2[q,q,q^4,q^6;q^{12}]_\infty \\ 
& = \frac{E^2(q^6) E^6(q^2)}{E^5(q^{12}) E(q^4) E^2(q)}.
\end{split}
\mylabel{eq:4.30} 
\eeq
We can simplify \eqn{4.30} with the aid of \eqn{4.17} as
\beq
[q^4,q^6;q^{12}]_\infty ([q^5;q^{12}]_\infty+q[q;q^{12}]_\infty)^2 = \frac{E^2(q^6) E^6(q^2)}{E^5(q^{12}) E(q^4) E^2(q)}.
\mylabel{eq:4.31} 
\eeq
Next, we use \eqn{4.18} to reduce \eqn{4.31} to the following easily verifiable identity
\beq
\frac{[q^2;q^{12}]_\infty^4[q^4,q^6;q^{12}]_\infty^3}{[q,q^3,q^5;q^{12}]_\infty^2} 
= \frac{E^2(q^6) E^6(q^2)}{E^5(q^{12}) E(q^4) E^2(q)}.
\mylabel{eq:4.32} 
\eeq
This completes our proof of \eqn{4.22}, \eqn{4.20}. 
We have \eqn{4.9}, as desired.

The proof of \eqn{4.10} is analogous. Again, we view this identity as a constant term in $z$ of the following
\beq
\sum_{j=14}^{19}s_j(z,q)=R_2(q)\sum_{n=-\infty}^\infty q^{9\frac{n(n+1)}{2}}z^n,
\mylabel{eq:4.33} 
\eeq
where
\beq
R_2(q)=q\frac{E^2(q^9) E^4(q^{12})}{E(q^3) E(q^6) E(q^{18})}.
\mylabel{eq:4.34} 
\eeq
Again, \eqn{4.24} holds true for the following $(i,j)$ pairs:
$(14,15)$, $(15,16)$, $(16,17)$, $(17,14)$, $(18,19)$, $(19,18)$.\\
And so, both sides of \eqn{4.33} satisfy \eqn{4.26}. Again it remains to show that \eqn{4.33} holds at one nontrivial point, say $z_1=-q^{-3}$. Observing that 
\beq
s_{14}(z_1,q)=s_{15}(z_1,q)=s_{19}(z_1,q)=0,
\mylabel{eq:4.35} 
\eeq
we find that \eqn{4.33} with $z=z_1$ becomes  
\beq
s_{16}(z_1,q)+s_{17}(z_1,q)+s_{18}(z_1,q)=R_2(q)\sum_{n=-\infty}^\infty q^{9\frac{n(n+1)}{2}}z_1^n,
\mylabel{eq:4.36} 
\eeq
Again, making repeated use of \eqn{4.15} and replacing $q^3$ by $q$, we can rewrite \eqn{4.36} as
\beq
\begin{split}
[q^3,q^4,q^6;q^{12}]_\infty([q^5;q^{12}]_\infty & - q[q;q^{12}]_\infty)+ q[q^2,q^3,q^3,q^4;q^{12}]_\infty \\
& = \frac{E^4(q^4) E^2(q^3)}{E^4(q^{12}) E(q^6) E(q^2)}.
\end{split}
\mylabel{eq:4.37} 
\eeq
If we multiply both sides of \eqn{4.37} by $\frac{[q^2;q^{12}]_\infty}{[q^4;q^{12}]_\infty}$ and take
advantage of \eqn{4.17} we find that
\beq
[q,q^5,q^6,q^6;q^{12}]_\infty + q[q^2,q^3;q^{12}]_\infty^2 = \frac{E^4(q^4) E^2(q^3)}{E^4(q^{12}) E(q^6) E(q^2)}
\frac{[q^2;q^{12}]_\infty}{[q^4;q^{12}]_\infty},
\mylabel{eq:4.38} 
\eeq
which is easy to recognize as \eqn{4.19}. This completes our proof of \eqn{4.33} and \eqn{4.10}.\\
To prove \eqn{4.11}, \eqn{4.12} we will follow well trodden path and observe that these identities 
are just constant terms in $z$ of
\beq
\sum_{j=20+4\alpha}^{23+4\alpha}s_j(z,q)=q^2\frac{E^2(q^9) E(q^{12}) E(q^{36})}{E(q^3)}\cdot
\sum_{n=-\infty}^\infty q^{9\frac{n(n+1)}{2}}z^n,
\mylabel{eq:4.39} 
\eeq
with $\alpha=0$ and $1$, respectively.\\
To prove that both sides of \eqn{4.39} satisfy \eqn{4.26} we verify that \eqn{4.24} 
holds for the following $(i,j)$ pairs:
$(20+4\alpha,21+4\alpha)$, $(21+4\alpha,22+4\alpha)$, $(22+4\alpha,23+4\alpha)$, $(23+4\alpha,20+4\alpha)$ with
$\alpha=0,1$. It remains to verify \eqn{4.39} at
\beqs
\tilde z_\alpha=-q^{6(1-2\alpha)},\quad\quad\quad \alpha=0,1.
\eeqs
Taking into account that
\beqs
s_{j+4\alpha}(\tilde z_\alpha,q)=0 
\eeqs
for $j=20,21,22$ and $\alpha=0,1$, we find that
\beqs
s_{23+4\alpha}(\tilde z_\alpha,q)=(-1)^{\alpha+1}q^{4+6\alpha}E^2(q^9)E(q^{12})E(q^{36}), 
\eeqs
which is easy to prove with the aid of \eqn{4.15}.\\
This completes our proof of \eqn{4.11} and \eqn{4.12}.

\bigskip

\section{Concluding Remarks}\label{sec:5}

\medskip

Making use of the Littlewood decomposition of $\pitc$
into its $s$-core and $s$-quotient,
\beqs
\tilde\phi(\pitc)=(\pi_{s\mbox{-core}},\hat\pi_0,\hat\pi_1,\ldots,\hat\pi_{s-1}),
\eeqs
together with 
\beqs
1+\omega_s+\omega_s^2+\cdots+\omega_s^{s-1}=0,
\eeqs
it is not hard to see that 
\beqs
\mgbgr(\pitc,s)=\mgbgr(\pi_{s\mbox{-core}},s).
\eeqs
In a recent paper \cite{O}, Olsson proved a somewhat unexpected result:
\begin{theorem}\label{t4}
Let $s,t$ be relatively prime positive integers, then the $s$-core of 
a $t$-core is, again, a $t$-core.
\end{theorem}
\noindent
In \cite{A}, Anderson established
\begin{theorem}\label{t5}
Let $s,t$ be relatively prime positive integers, then the number of partitions, which are simultaneously $s$-core and $t$-core is $\frac{\bin{s+t}{s}}{s+t}$.
\end{theorem}
\noindent
Remarkably, the three observations above imply our Theorem 1.2. \\
Moreover, our Theorem 1.3 implies 
\begin{corollary}
Let $s,t$ be relatively prime positive integers. 
Then no two distinct $(s,t)$-cores share the same value of $\mgbgr\mymod{s}$, 
when $s$ is prime, or when $s$ is composite and $t<2p_s$, where $p_s$ is a smallest prime divisor of $s$.
\end{corollary}
On the other hand, when the conditions on $s$ and $t$ in the corollary above are not met, two distinct
$(s,t)$-cores may, in fact, share the same value of $\mgbgr\mymod{s}$. For example, consider two relatively 
prime integers $s$ and $t$ such that $2\mid s$, $s>2$, $t>1+\frac{s}{2}$, $t\neq s+1$.
In this case partition $[1^{\frac{s}{2}-1},2,1+\frac{s}{2}]$ and empty partition $[\;\;]$ 
are two distinct $(s,t)$-cores such that
\beqs
\mgbgr\big(\big[1^{\frac{s}{2}-1},2,1+\frac{s}{2}\big],s\big)=\mgbgr([\;\;],s)=0.
\eeqs

\bigskip
\noindent
\textbf{Acknowledgement}
\medskip

\noindent
We are grateful to Hendrik Lenstra for his contribution to the proof of Lemma 1.4.
We would like to thank Christian Krattenthaler, Jorn Olsson and Peter Paule for their kind interest and helpful discussions.

\bigskip
\bigskip

\bibliographystyle{amsplain}

\end{document}